# HARMONIC MOTION AND CASSINI OVALS

March 23, 2013

**Khristo Boyadzhiev**

Department of Mathematics and Statistics

Ohio Northern University

Ada, OH 45810

k-boyadzhiev@onu.edu

**Irina Boyadzhiev**

Department of Mathematics

Ohio State University (Lima)

Lima, OH 45804

boyadzhiev.1@osu.edu

**Abstract** We consider a two-dimensional free harmonic oscillator where the initial position is fixed and the initial velocity can change direction. All possible orbits are ellipses and their enveloping curve is an ellipse too. We show that the locus of the foci of all elliptical orbits is a Cassini oval. Depending on the magnitude of the initial velocity we observe all three kinds of Cassini ovals, one of which is the lemniscate of Bernoulli. These Cassini ovals have the same foci as the enveloping ellipse.

## 0. INTRODUCTION

The main result in this paper is about two-dimensional harmonic oscillators. The trajectories of the oscillating points are ellipses depending on a parameter. In Section 3 we prove that the locus of the foci of these ellipses is a Cassini oval. The form of this oval depends on the magnitude of the initial velocity. We approach this theorem gradually, by discussing first enveloping curves to such periodic motions. The material in the first two sections is not new and these sections can be considered a review where important details are discussed thoroughly. Section 1 has the purpose to introduce the reader into the specific subject matter. Also, at the end of the first section we point out one interesting fact - the locus of the foci of the set of parabolas is a circle.

## 1. PARABOLA OF SAFETY

Suppose a projectile is thrown from the center of the coordinate system $O(0,0)$ with initial velocity $v_0$ at angle $\alpha$ with the $x$-axis. The motion of the projectile in vacuum is described by the parametric equations

$$x = (v_0 \cos \alpha)t, \quad y = (v_0 \sin \alpha)t - \frac{g}{2}t^2 \tag{1.1}$$

($g$ - the gravitational constant). The trajectory is a parabola with Cartesian equation

$$y = x \tan \alpha - \frac{g}{2v_0^2 \cos^2 \alpha} x^2 \tag{1.2}$$

Setting $y = 0$ and solving for $x$ we find the second $x$-intercept, $x = \frac{v_0^2}{g} \sin(2\alpha)$. This shows that the maximal reachable distance on the positive $x$-axis is $x_{max} = \frac{v_0^2}{g}$ corresponding to

$\alpha = \frac{\pi}{4}$. Symmetrically, the point $(-x_{max}, 0)$ will be reached when $\alpha = \frac{3\pi}{4}$.

When the angle $\alpha$ changes from $0$ to $\pi$ we obtain a family of parabolas with the same initial velocity. They cover a certain region in the $xy$-plane.

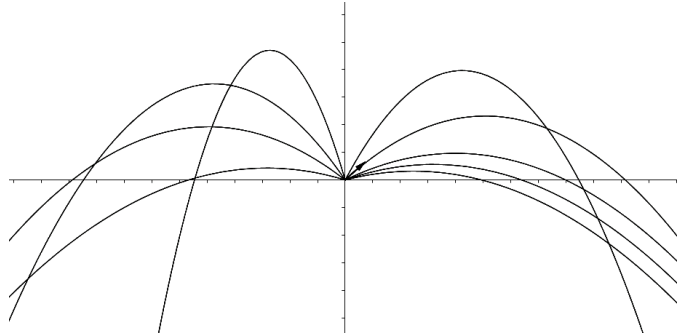

**Figure 1 Some trajectories**

The question is: What is the boundary of this region? What is the "safety" curve beyond which all points are "safe"?

We shall find the equation of this curve and see that it is a parabola. This parabola is called "parabola of safety". Notice that the points $(\pm x_{max}, 0)$ should be on this parabola.

First we rewrite equation (1.2) as a quadratic equation in $\tan \alpha$ by using that $\frac{1}{\cos^2 \alpha} = 1 + \tan^2 \alpha$

$$\frac{gx^2}{2v_0^2} \tan^2 \alpha - x \tan \alpha + y + \frac{gx^2}{2v_0^2} = 0. \tag{1.3}$$

Next we reason this way: If a point $(x, y)$ can be reached by two parabolas of the above family, then it is not on the "safety" curve. Therefore, if $(x, y)$ is on the "safety" curve, the equation (1.3) has only one solution for $\tan\alpha$. In that case the discriminant of equation (1.3) is zero, i.e.

$$x^2 - \frac{2gx^2}{v_0^2}(y + \frac{gx^2}{2v_0^2}) = 0,$$

which we can reduce by $x^2$ and solve for $y$

$$y = \frac{-g}{2v_0^2}x^2 + \frac{v_0^2}{2g}. \tag{1.4}$$

This is a parabola with vertex $V = (0, \frac{v_0^2}{2g})$, focus at the origin $O(0,0)$, and $x$-intercepts

$$x = \pm\frac{v_0^2}{g} = \pm x_{max}.$$

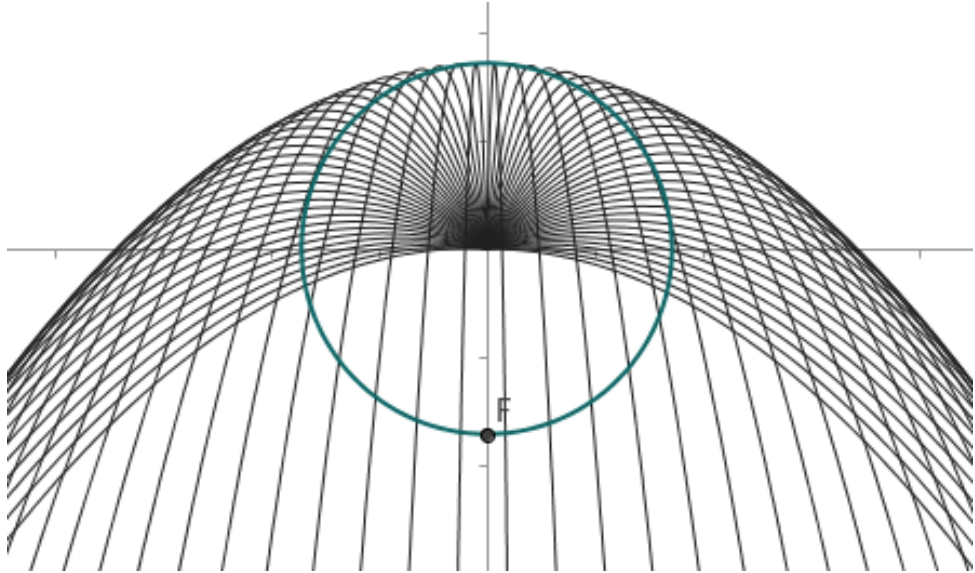

**Figure 2 Enveloping parabola of safety and the locus of the foci**

The fact that the "safety" curve is a parabola is not new; it was discovered in the 17[th] century by the Italian mathematician Evangelista Torricelli by purely geometrical means.

The parabola (1.4) touches every trajectory (1.2) (or (1.3)) only once, i.e. it is the enveloping curve for all these parabolas. Replacing $y$ from (1.4) in (1.3) and solving for $\tan\alpha$ we find the equation

$$\tan \alpha = \frac{v_0^2}{g\,x}, \tag{1.5}$$

connecting the angle $\alpha$ and the $x$-coordinate of the intercept. Here $x = 0$ corresponds to $\alpha = \frac{\pi}{2}$ and for every $x \neq 0$ there is only one angle $\alpha \neq \frac{\pi}{2}$ in the interval $0 < \alpha < \pi$ satisfying the above equation. The endpoints $\alpha = 0$ and $\alpha = \pi$ correspond to $x = +\infty$ and $x = -\infty$.

The parabola

$$y = -\frac{g}{2v_0^2} x^2 \tag{1.6}$$

corresponding to $\alpha = 0$, (the same parabola also results from $\alpha = \pi$) is asymptotic to the enveloping parabola (1.4) at $x = \pm\infty$. All trajectories under it (when $-\pi < \alpha < 0$ or $\pi < \alpha < 2\pi$) do not touch the enveloping parabola.

Rotating the parabola of safety (1.4) about the $y$-axis generates a paraboloid with the same vertex and same focus. This is the paraboloid of safety. The points outside this paraboloid can not be reached by projectiles from the origin with the same initial velocity in 3-dimensions.

**Remark 1**. Another interesting geometric object is the locus of the foci of the parabolas (1.2). For each parabola with equation (2) the focus has coordinates

$$F\left(\frac{v_0^2}{2g}\sin 2\alpha,\ -\frac{v_0^2}{2g}\cos 2\alpha\right). \tag{1.7}$$

The trajectory of these point when $\alpha$ changes from $0$ to $\pi$ is a circle with center at the origin and radius $R = \frac{v_0^2}{2g}$ see Figure 1. In three dimensions we have a sphere (inside the paraboloid of safety) generated by rotating this circle about the $y$-axis.

References for this section – [1], [2], [3], [6], [8].

## 2. ELLIPSE OF SAFETY IN HARMONIC MOTION

Recall that the free harmonic oscillator is governed by the differential equation

$$x'' + \omega^2 x = 0 \tag{2.1}$$

where $x(t)$ is the position function and $t \geq 0$ is time. This equation results from Newton's law $F = ma$ by using the force $F = -kx$, where $k > 0$ is a constant. Then we have $mx'' = -kx$ or (1) with $\omega^2 = k/m$.

Suppose the initial position is the point $(x_0, 0)$ in the $xy$-plane and the initial velocity is $v_0 = x'(0)$. The law of motion is

$$x = x_0 \cos(\omega t) + \frac{v_0}{\omega} \sin(\omega t). \tag{2.2}$$

The point $M(x(t), 0)$ oscillates over the interval $[-c, c]$, where $c = \sqrt{x_0^2 + \frac{v_0^2}{\omega^2}}$.

Now we extend this motion to two dimensions by allowing the initial velocity to be a vector $\bar{v}_0 = (v_0 \cos\alpha, v_0 \sin\alpha)$, cutting angle $\alpha$ with the $x$-axis. The motion is described by the vector function $\bar{r}(t) = (x(t), y(t))$ whose coordinates satisfy the equations.

$$\begin{aligned} x'' + \omega^2 x &= 0 \\ y'' + \omega^2 y &= 0 \end{aligned} \tag{2.3}$$

with initial position $(x_0, 0)$ and initial velocity $\bar{v}_0 = (v_0 \cos\alpha, v_0 \sin\alpha)$. Solving for $x$ and $y$ separately we find the parametric equations

$$x = x_0 \cos(\omega t) + \frac{v_0}{\omega} \cos\alpha \sin(\omega t), \quad y = \frac{v_0}{\omega} \sin\alpha \sin(\omega t), \tag{2.4}$$

and the trajectory is an ellipse. Setting for convenience $p = \frac{v_0}{\omega}$ and replacing in the first equation

$$\sin(\omega t) = \frac{y}{p \sin\alpha} \tag{2.5}$$

we find $x = x_0 \cos(\omega t) + y \cot\alpha$ and then we solve for $\cos(\omega t)$,

$$\cos(\omega t) = \frac{x - y \cot\alpha}{x_0}. \tag{2.6}$$

Now from (2.5) and (2.6) we find the equation of the ellipse in Cartesian coordinates,

$$\left(\frac{x - y \cot\alpha}{x_0}\right)^2 + \left(\frac{y}{p \sin\alpha}\right)^2 = 1 \tag{2.7}$$

This is an ellipse with center $(0,0)$ (notice the symmetry according to $(0,0)$). One simple representative of this family is the ellipse corresponding to $\alpha = \dfrac{\pi}{2}$,

$$\frac{x^2}{x_0^2} + \frac{y^2}{p^2} = 1,$$

with radii $x_0$ and $p = v_0/\omega$.

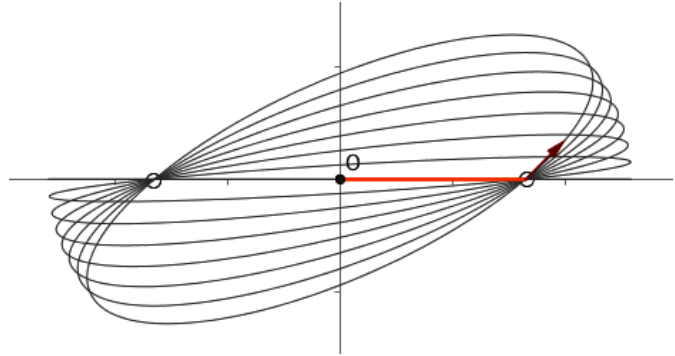

**Figure 3 Several orbits with the same speed, $v_0 < x_0$**

Now we shall try to identify the region $G$ filled by all these ellipses when we change $\alpha$ in $(0, \pi)$. For this purpose we rewrite equation (2.7) as a quadratic equation for $\cot \alpha$.

With the help of the identity $\dfrac{1}{\sin^2 \alpha} = 1 + \cot^2 \alpha$ the equation becomes

$$(x - y \cot \alpha)^2 p^2 + x_0^2 y^2 (1 + \cot^2 \alpha) - x_0^2 p^2 = 0, \text{ or}$$

$$y^2(p^2 + x_0^2)\cot^2 \alpha - 2xyp^2 \cot \alpha + (x^2 - x_0^2)p^2 + x_0^2 y^2 = 0. \tag{2.8}$$

We reason that a point $(x, y)$ is on the boundary of $G$ if it can be reached by only one ellipse in the above set of ellipses. This means that for points $(x, y)$ on the boundary, equation (2.8) has only one solution for $\cot \alpha$. In this case the discriminant is zero, i.e.

$$4x^2 y^2 p^4 - 4y^2 (p^2 + x_0^2)[(x^2 - x_0^2)p^2 + x_0^2 y^2] = 0.$$

Simplifying this we obtain

$$\frac{x^2}{x_0^2 + p^2} + \frac{y^2}{p^2} = 1 \tag{2.9}$$

which is an ellipse with center at the origin, foci at $(-x_0, 0)$ and $(x_0, 0)$, big axis $\sqrt{x_0^2 + p^2}$ and small axis $p$, i.e. the axes are $a = \sqrt{x_0^2 + \dfrac{v_0^2}{\omega^2}} = \sqrt{x_0^2 + \dfrac{mv_0^2}{k}}$ (which is the maximal amplitude) and $p = \dfrac{v_0}{\omega} = v_0 \sqrt{\dfrac{m}{k}}$. The enveloping ellipse can be seen in figures 4, 5, and 6.

If we rotate this ellipse about the $x$-axis, it will generate the ellipsoid of safety, with the same center and same foci. Beyond this ellipsoid is the "safe" space, where the oscillating point $M$ with initial speed $|\bar{v}_0|$ cannot reach.

For the ellipse of safety and other similar results see [5].

### 3. THE TRAJECTORY OF THE FOCI

When $\alpha$ varies the foci of the orbits (2.4) trace a remarkable curve, a Cassini oval. The Cassini oval is defined as the locus of all points $(x, y)$ whose distances to two fixed points (foci) $(-\lambda, 0)$ and $(\lambda, 0)$ have a constant product $\mu^2$, i.e.

$$[(x+\lambda)^2 + y^2][(x-\lambda)^2 + y^2] = \mu^4$$

We have the following theorem where without loss of generality we assume that the frequency $\omega = 1$.

**Theorem**. *The foci of the elliptical trajectories*

$$x = x_0 \cos(t) + v_0 \cos\alpha \sin(t), \quad y = v_0 \sin\alpha \sin(t), \quad 0 \leq t \leq 2\pi,$$

*with initial position $(x_0, 0)$ and initial velocity $\bar{v}_0 = (v_0 \cos\alpha, v_0 \sin\alpha)$ when $\alpha$ changes from $\alpha = -\pi$ to $\alpha = \pi$ trace a Cassini oval with Cartesian equation*

$$[(x+x_0)^2 + y^2][(x-x_0)^2 + y^2] = v_0^4 \tag{3.1}$$

*or, equivalently,*

$$(x^2 + y^2)^2 - 2x_0^2(x^2 - y^2) = v_0^4 - x_0^4 \tag{3.2}$$

*and polar equation*

$$r^4 - 2x_0^2 r^2 \cos 2\theta = v_0^4 - x_0^4. \tag{3.3}$$

*This Cassini oval sits symmetrically inside the ellipse of safety. In particular, when $x_0 = v_0$ the locus of the foci is Bernoulli's lemniscate with Cartesian equation*

$$(x^2 + y^2)^2 = 2x_0^2(x^2 - y^2) \tag{3.4}$$

*and polar equation*

$$r^2 = 2x_0^2 \cos 2\theta, \quad -\frac{\pi}{4} \leq \theta \leq \frac{\pi}{4}. \tag{3.5}$$

*Moreover, in this case $\alpha = 2\theta$, where $\theta$ is the polar angle. Thus equation (3.4) becomes*

$$r^2 = 2x_0^2 \cos\alpha, \quad -\frac{\pi}{2} \leq \alpha \leq \frac{\pi}{2}. \tag{3.6}$$

**Remark 2.** Although in equation (3.6) $\cos\alpha$ is only non-negative, allowing the polar radius in this equation to be negative we obtain the entire lemniscate.

**Remark 3.** As stated in the theorem, when $v_0 = x_0$ the Cassini oval is a lemniscate. The other two cases $v_0 < x_0$ and $v_0 > x_0$ provide two-pieces and one-piece ovals correspondingly. All three cases of Cassini ovals are presented in figures 4, 5, and 6 below. In these figures $A = (-x_0, 0)$ and $B = (x_0, 0)$ are the foci of the enveloping ellipse and also the foci of the Cassini ovals. The points $F_1$ and $F_2$ are the foci and also the endpoints of the degenerated ellipse-segment $[-a, a]$, $a = \sqrt{x_0^2 + v_0^2}$ occurring for $\alpha = 0$. These are the only two common points for the Cassini ovals and the enveloping ellipse.

Notice that equation (3.1) is a bipolar equation expressing the fact that if $F(x, y)$ is a point on the Cassini oval, the product of its distances to $A$ and $B$ is $v_0^2$.

Good reference for Cassini ovals are [4], [7], [9].

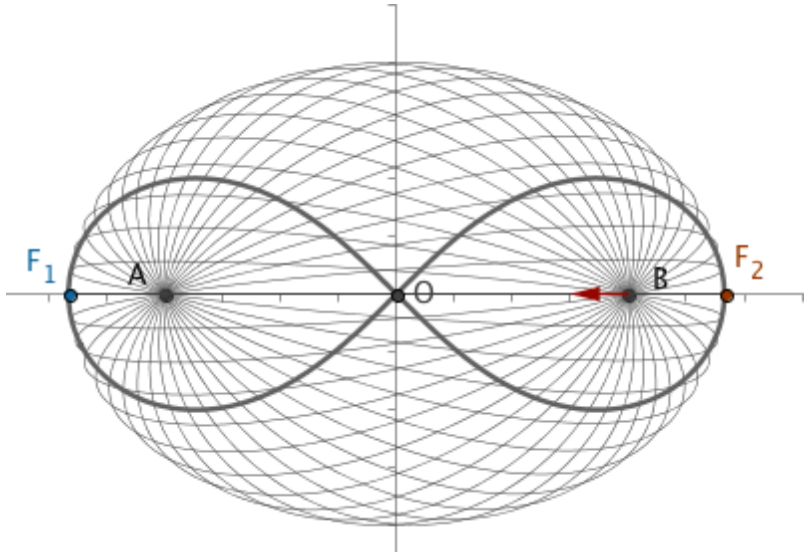

**Figure 4 Lemniscate of Bernoulli,** $v_0 = x_0$

When $v_0 < x_0$ the Cassini curve consists of two ovals, as shown on Figure 5.

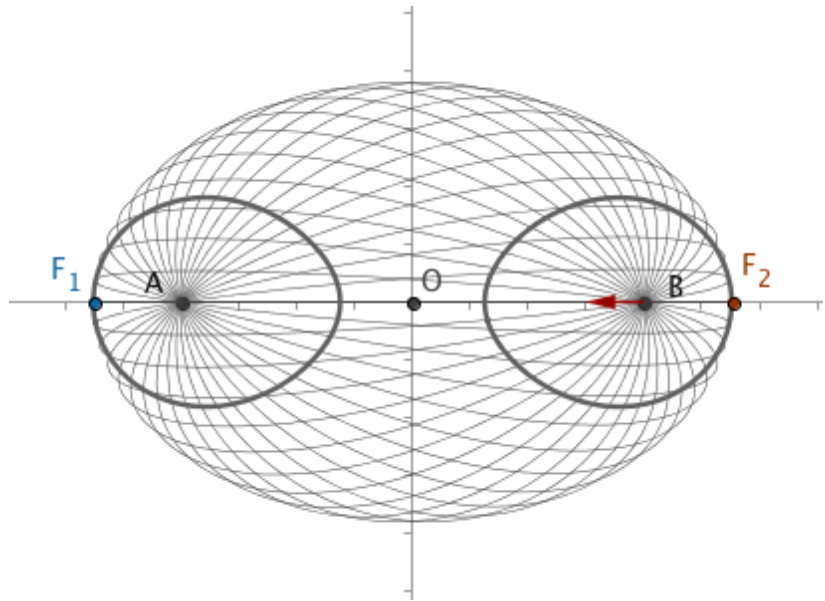

**Figure 5 Two Cassini ovals,** $v_0 < x_0$

When $v_0 > x_0$ the Cassini oval consists of one piece.

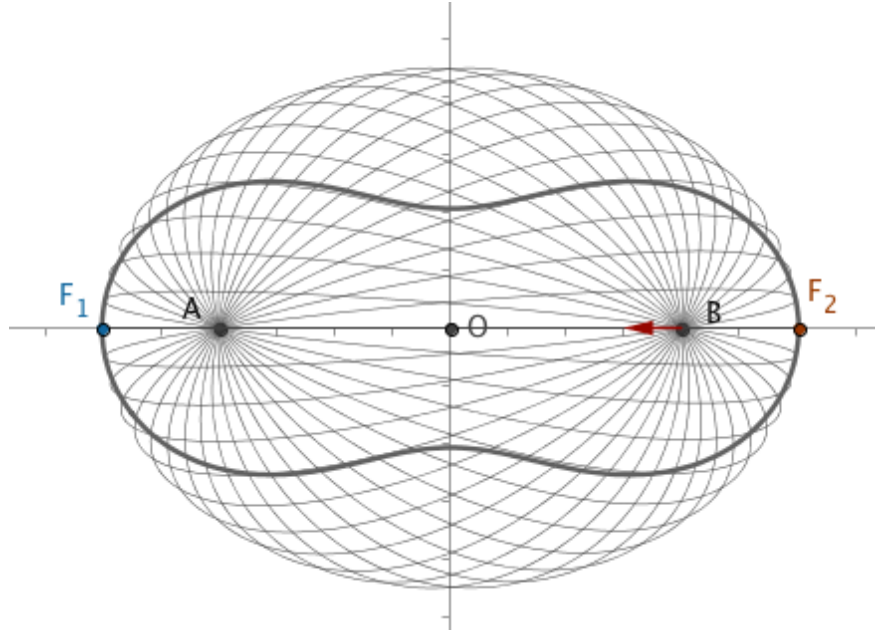

**Figure 6 Cassini oval** $v_0 > x_0$

*Proof of the theorem.* We first give a proof for the most interesting case $x_0 = v_0$. For this case it is easy to write explicit parametric equations for the trajectory of the foci. The proof of the general case will be different.

When $x_0 = v_0$ and $\omega = 1$ equations (2.4) become

$$x = x_0(\cos t + \cos \alpha \sin t), \quad y = x_0 \sin \alpha \sin t .  \tag{3.7}$$

From here we compute

$$r^2 = x^2 + y^2 = x_0^2 (\cos^2 t + \sin^2 t \cos^2 \alpha + \sin^2 t \sin^2 \alpha + \cos \alpha \sin 2t)$$

That is,

$$r^2 = x_0^2 (1 + \cos \alpha \sin 2t) \tag{3.8}$$

Assume first that $-\dfrac{\pi}{2} \le \alpha \le \dfrac{\pi}{2}$, so that $\cos \alpha \ge 0$. The maximal and minimal values of $r^2$ are the squares of the big axis $a$ and the small axis $b$ of the ellipse (3.7). They are obtained when $\sin 2t = \pm 1$ $(t = \dfrac{\pi}{4}, \dfrac{3\pi}{4}, ...)$. Thus

$$a^2 = x_0^2 (1 + \cos \alpha), \quad b^2 = x_0^2 (1 - \cos \alpha) .  \tag{3.9}$$

For the focal distance $c$ we have $c^2 = a^2 - b^2$, so that

$$c^2 = 2x_0^2 \cos\alpha \tag{3.10}$$

The vertex of the ellipse in the first quadrant happens when $t = \dfrac{\pi}{4}$ and the coordinates of this vertex are (from (3.7).

$$V\left(\frac{x_0(1+\cos\alpha)}{\sqrt{2}}, \frac{x_0 \sin\alpha}{\sqrt{2}}\right). \tag{3.11}$$

Multiplying by $c/a$ these coordinates we find the coordinates of the focus in the first quadrant as functions of $\alpha$,

$$x = x_0\sqrt{\cos\alpha(1+\cos\alpha)}, \quad y = \frac{x_0 \sin\alpha\sqrt{\cos\alpha}}{\sqrt{1+\cos\alpha}}. \tag{3.12}$$

From here $x^2 + y^2 = 2x_0^2 \cos\alpha$ confirming (3.10) and also

$$x^2 - y^2 = 2x_0^2 \cos^2\alpha. \tag{3.13}$$

Equation (3.4) follows immediately from here. The restriction $\cos\alpha \geq 0$ is not essential, because equation (3.4) extends by symmetry to all quadrants.

Let now $\theta$ be the polar angle and $x = r\cos\theta$, $y = r\sin\theta$ the standard polar relations. Then

$$x^2 - y^2 = r^2(\cos^2\theta - \sin^2\theta) = r^2 \cos 2\theta$$

and comparing this to (3.13) and (3.10) for the coordinates of the focus we find that $\cos\alpha = \cos 2\theta$ and thus $\alpha = 2\theta$. The proof is completed.

**Remark 4.** Notice that the vertex $V$ in (3.11) traces a semicircle for $-\dfrac{\pi}{2} \leq \alpha \leq \dfrac{\pi}{2}$ with radius $\dfrac{x_0}{\sqrt{2}}$ and center $\left(\dfrac{x_0}{\sqrt{2}}, 0\right)$. The vertex cannot complete the whole circle when $\alpha$ moves beyond $\pm\dfrac{\pi}{2}$, because at $\alpha = \pm\dfrac{\pi}{2}$ the ellipse (3.4) becomes a circle with $a = b$ and the next moment the two axes $a$ and $b$ change places. The second half of the circle is centered at $\left(-\dfrac{x_0}{\sqrt{2}}, 0\right)$, see figure 7.

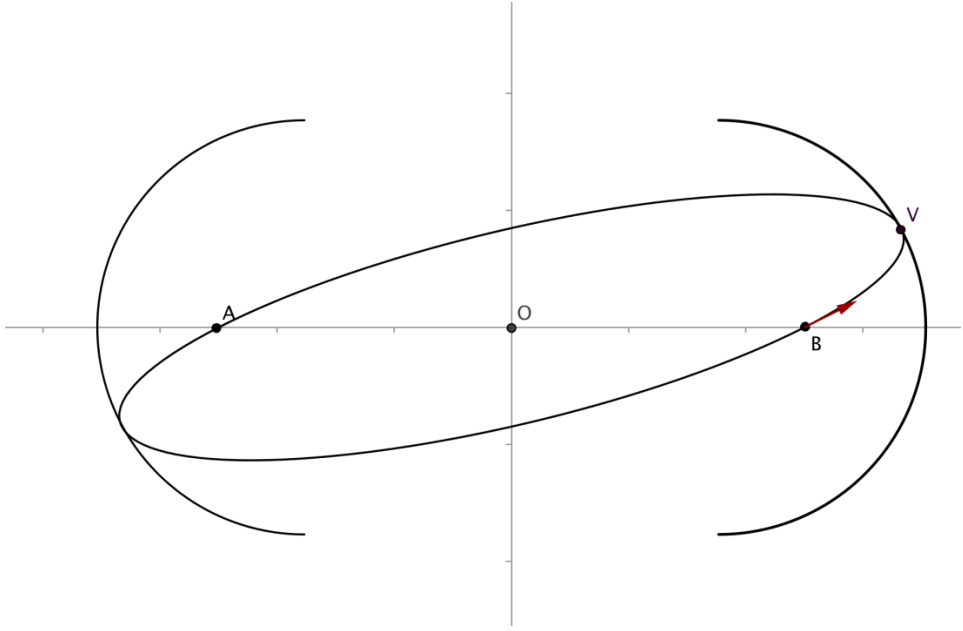

**Figure 7 Trajectories of the vertices**

*Proof of the theorem in the case $v_0 \neq x_0$.*

We need to prove equation (3.1). When $v_0 \neq x_0$ we compute from (2.4)

$$r^2 = x^2 + y^2 = x_0^2 \cos^2 t + v_0 \sin^2 t + 2x_0 v_0 \cos\alpha \sin t \cos t$$

which can be written as

$$r^2 = \frac{1}{2}(x_0^2 + v_0^2) + \frac{1}{2}(x_0^2 - v_0^2)\cos 2t + v_0 \sin^2 t + x_0 v_0 \cos\alpha \sin 2t \ . \tag{3.14}$$

Maximal and minimal values occur when

$$\tan 2t = \frac{2x_0 v_0 \cos\alpha}{x_0^2 - v_0^2} \ .$$

If maximum occurs for some $t$, then the next extremum, the minimum, will happen for $t + \frac{\pi}{2}$ as $\tan 2t$ is the same. Thus for the two radii $a$ and $b$ of the ellipse (2.4) we have

$$a^2 = \frac{1}{2}(x_0^2 + v_0^2) + \frac{1}{2}(x_0^2 - v_0^2)\cos 2t + x_0 v_0 \cos\alpha \sin 2t$$

$$b^2 = \frac{1}{2}(x_0^2 + v_0^2) - \frac{1}{2}(x_0^2 - v_0^2)\cos 2t - x_0 v_0 \cos\alpha \sin 2t$$

and from here we obtain the remarkable equation

$$a^2 + b^2 = x_0^2 + v_0^2 \ . \tag{3.15}$$

We use now the construction on Figure 7 with an arbitrary elliptical orbit.

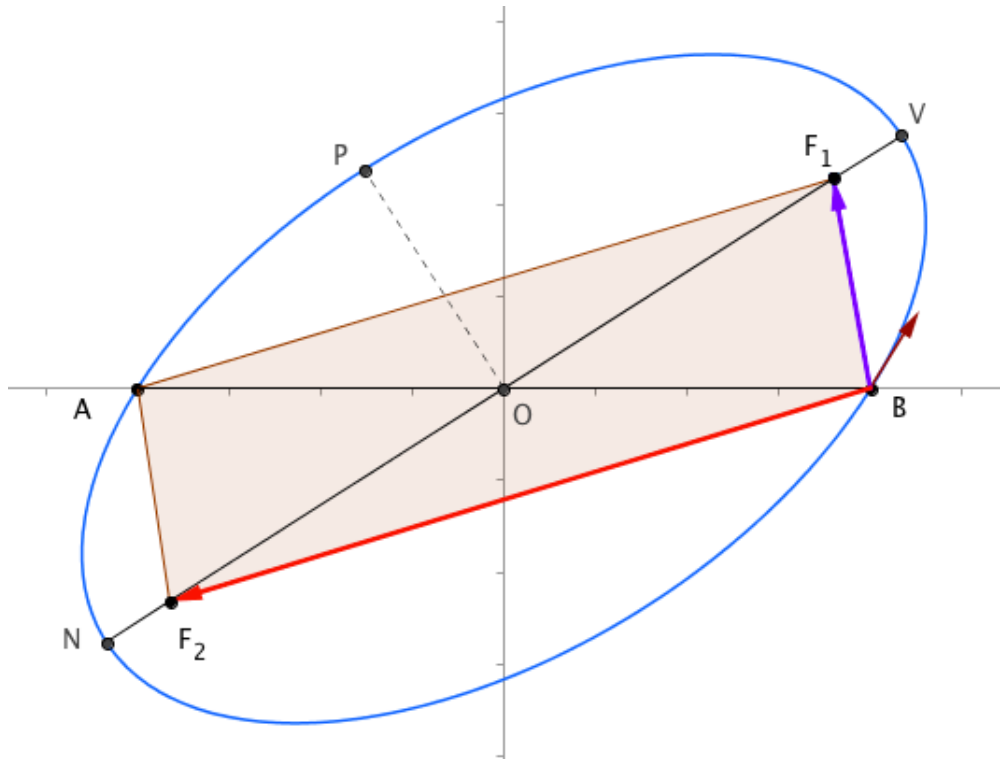

**Figure 8 An important parallelogram**

As before, $A = (-x_0, 0)$, $B = (x_0, 0)$. The quadrilateral $AF_2 BF_1$ is a parallelogram – note that the diagonals $AB$ and $F_1 F_2$ cut each other in half. Thus $AF_2 = BF_1$ and $AF_1 = BF_2$. From the parallelogram law we have

$$AB^2 + F_1 F_2^{\,2} = 2(F_1 B^2 + F_2 B^2)$$

Note that $AB = 2x_0$ and $F_1 F_2^{\,2} = 4(a^2 - b^2)$, so that

$$F_1B^2 + F_2B^2 = 2OF_1^2 + 2x_0^2 \tag{3.16}$$

At the same time, using the property of the ellipse,

$$4a^2 = (F_1A + F_1B)^2 = F_1A^2 + F_1B^2 + 2F_1AF_1B$$

Therefore, with the help of (3.16) and using that $F_1A = F_2B$,

$$2F_1AF_1B = 4a^2 - F_2B^2 - F_1B^2 = 4a^2 - 2(a^2 - b^2) - 2x_0^2 = 2(a^2 + b^2) - 2x_0^2$$

and now using (3.15),

$$F_1AF_1B = v_0^2$$

which is equation (3.1). Done!

**Remark 3**. On Figure 7, $a = OV$ and $b = OP$. Equation (3.15) shows that the segment $PV$ has constant length $\sqrt{x_0^2 + v_0^2}$ for all elliptical orbits. This quantity is exactly the large radius of the enveloping ellipse (2.9) when $\omega = 1$. Moreover, the enveloping ellipse and all Cassini ovals have the same foci $A$ and $B$.

**References**

[1]     **Denis Donnely,** The parabolic envelope of constant initial speed trajectories, *Am. J. Phys*. 60(12) (1992), 1149-1150.

[2]     **Richard Fitzpatrick**, *Newtonian Dynamics*, Lecture Notes, 2011.

[3]     **Derek Hart, Tony Croft,** *Modeling with Projectiles,* Wiley, 1988.

[4]     **J. Dennis Lawrence**, *A Catalog of Special Plane Curves*, Dover 1972.

[5]     **Jean-Marc Richard,** Safe domain and elementary geometry, *European Journal of Physics*, 25 (2004), 835-848.

[6]     **John L. Synge, Byron A. Griffith**, *Principles of Mechanics*, McGraw-Hill, 1959.

[7]     **F. Gomes Teixeira**, *Traité des Courbes Spéciales Remarquables*, Tome 1, Chelsea, 1971.

]8]     **R.E. Warner and L.A. Huttar**, The parabolic shadow of a Coulomb scatter, *Am. J. Phys*. 59(8) (1991),755-756.

[9]     **Robert C. Yates**, *Curves and Their Properties*, 1952.